\documentclass[twoside,11pt,reqno]{amsart}
\usepackage{amsfonts,amsmath,amssymb,amscd,mathrsfs}

\usepackage{amsfonts, amsbsy, amsmath, amsthm, amssymb, latexsym, verbatim, enumerate, color}
\definecolor{re}{rgb}{1,0.2,0.2}           
 \definecolor{gr}{rgb}{0,1,0}
 \definecolor{bl}{rgb}{0,0,0.6}
  \definecolor{bl2}{rgb}{0,1,0}

\usepackage{youngtab}
\newtheorem{propo}{Proposition}[section]

\newtheorem{lemma}[propo]{Lemma}
\newtheorem{corol}[propo]{Corollary}
\newtheorem{theor}[propo]{Theorem}
\newtheorem{theo}[propo]{Theorem}

\newcommand{\ld}{,\ldots ,}

\newcommand{\lan}{ \langle }
\newcommand{\ran}{ \rangle }

\newcommand{\Irr}{\mathop{\rm Irr}\nolimits}

\newcommand{\al}{\alpha}

\newcommand{\dne}{\hfill $\Box$ \vspace{0.3cm}}

\newcommand{\up}{^{-1}}

\newcommand{\supp}{{\rm supp} }

\def\d12{{_{12}}}

\def\ii{{if and only if }}

\def\ir{{irreducible }}

\def\irr{{irreducible representation }}

\def\itf{{It follows that }}

\def\rep{{representation }}

\def\syl{{Sylow $p$-subgroup }}

\def\di{diagram }

\newcommand{\pf}{{\it Proof:\quad}}

\newcommand{\lp}{{\big( }}
\newcommand{\rp}{{\big) }}

\newcommand{\gaa}{\alpha}

\newcommand{\Ga}{\Gamma}

\newcommand{\gl}{\lambda}

\newcommand{\gO}{\Omega}

\newcommand{\re}{\mathbb{R}}  
  
\newcommand{\Fi}{\mathbb{F}}

\newcommand{\el}{\end{lemma}}
\newcommand{\om}{\omega }

\newcommand{\bl}{\begin{lemma}\label}
\date{}
\begin{document}

\title{Remarks  on  singular Cayley graphs and \\vanishing elements of simple groups}
\author{J. Siemons and A. Zalesski}
\address{A. Zalesski: Department of Physics, Mathematics and Informatics \\ National Academy of Sciences of Belarus 
\\Nezavisimosti prospekt 66\\
 220072 Minsk\\
  Belarus}
\email{
    alexandre.zalesski@gmail.com }

\address{J. Siemons: School of Mathematics, University of East Anglia, Norwich, NR4 7TJ, UK}
\email{j.siemons@uea.ac.uk}
\maketitle

{\sc Abstract:}\,   Let $\Ga$ be a finite graph 
 and let $A(\Gamma)$ be its 
 adjacency matrix.  Then $\Ga$ is {\it singular} if $A(\Ga)$ is  singular. The singularity  of graphs is of certain interest in  graph theory and algebraic combinatorics. Here we investigate this problem for Cayley graphs ${\rm Cay}(G,H)$ when $G$ is a finite group and when  the connecting set $H$ is a union of conjugacy classes of $G.$ 
In this situation the singularity problem reduces to finding an irreducible character $\chi$ of $G$ for which $\sum_{h\in H}\,\chi(h)=0.$
  At this stage we focus on the case when $H$ is a single conjugacy class $h^G$ of $G$; in this case the above equality is equivalent to $\chi(h)=0$. Much is known in this situation, with essential information coming from the block theory of representations of finite groups. An element $h\in G$ is called vanishing if $\chi(h)=0$ for some \ir character $\chi$ of $G.$ We study vanishing elements  mainly in finite simple groups and  in  alternating groups in particular. We  
suggest  some approaches for constructing singular Cayley graphs.
\footnote{{\sc Keywords:}\, Singular Cayley graphs, vertex transitive graphs, vanishing elements, block theory of symmetric and alternating groups. ~~~~\\
{\sc Mathematics Subject Classification} 05E99, 68R10, 20G30
}

\section{Introduction}
 Let $\Ga$ be a finite graph 
  and let $A(\Gamma)$ be its 
  adjacency matrix. 
 Then $\Ga$ is {\it singular} if $A(\Ga)$ is a singular matrix. Alternatively, $\Ga$ is singular if and only if its spectrum contains the eigenvalue $0.$ All graphs in this paper are undirected, without loops and without multiple edges; for all definitions please see Section~2.

Singular graphs play a significant role in  graph theory, and there are many applications in physics and chemistry, see Section~2. While the literature on graph spectra is vast it is not likely that a general theory of graph singularity {\it per se} will emerge. Some progress however can be made for  graphs which admit a group of  automorphisms that is transitive on the vertices of the graph.  In some cases the singularity problem then can be solved using techniques from ordinary character theory. The main purpose of this paper is to investigate these applications of character theory in graph theory. 


\medskip
In the following $G$ denotes a finite group 
and $H$ denotes a {\it connecting set} in $G.$ This is a subset of $G$ such that \,(i) $H$ does not contain the identity element $1$ of $G,$\,
(ii)\, $H=H^{-1}:=\{h^{-1}\,|\,h\in H\}$\, and\, (iii) $H$ generates $G$, that is, $H$ does not lie in any proper subgroup of $G$. 
From  these data the Cayley graph $\Ga={\rm Cay}(G,H)$ with vertex set $V=G$ and connecting sets $H$ can be defined, see Section 2. Here $\Ga$ is a regular graph of degree $|H|$ and the  group $G$ acts transitively on the vertices of $\Ga.$ Note though that a graph may be the Cayley graph of more than one group and connecting set. 

In this paper we specify the singularity problem to
Cayley graphs ${\rm Cay}(G,H)$ when 
the connecting set $H$ is 
$G$-invariant, that is,  $H$ is a union of conjugacy classes of $G.$ In this case the following theorem reduces the singularity problem to a problem of character
 theory. 

   \begin{theor}\label{SingCay-2}  Let $G$ be a group with a G-invariant connecting set $H$. 
    Then ${\rm Cay}(G,H)$ is singular if and only if there is an irreducible character $\chi$ of $G$ with $\sum_{h\in H}\,\chi(h)=0.$ In particular, if $H=h^{G}$ is a single conjugacy class then ${\rm Cay}(G,H)$ is singular if and only if there is an irreducible character $\chi$ of $G$ such that $\chi(h)=0.$ \end{theor}
 
{\sc Comments:\, }1. Burnside's theorem on character zeros \cite[\S 32, Exercise 3]{CR} shows that every character $\chi$ of degree $>1$ takes the value $\chi(h)=0$ for some $h\in G.$ Hence for every non-abelian group there exists a singular Cayley graph.  \\[5pt]
2. Obviously, if $\chi(h)=0$ for all $h\in H$ then ${\rm Cay}(G,H)$ is singular. 
 In addition, $\chi(h)=0$ \ii $\chi(h\up)=0$. Therefore, for constructing singular Cayley graphs it suffices to specify an \ir character $\chi$ of $G$  and a  set $X$ generating 
 $G$ such that $\chi$ takes the value $0$ on $X.$  Then, setting  $H=\cup_{ g\in G} \,\,g(X\cup X\up) g\up,$ we conclude that $H$ is a connecting set and  so ${\rm Cay}(G,H)$ is singular.  \\[5pt]
3. One may ask whether $\chi(h)+\chi(h\up)=0$ implies that $\chi(h)=0$. This is not so, see \cite{Atl}. In $G=PSU_3(3)$ there is an element $h\in G$ of order $4$ and two \ir characters $\chi$ of degree 28 so that $\chi(h)+\chi(h^{-1})=0$ while  $\chi(h)=\pm 4i\neq 0.$ 


 

 If the character table of a group $G$ is available explicitly (which is the case for  sporadic simple groups, say)  then one can determine  in principle all singular Cayley graphs  ${\rm Cay}(G,M\cup M^{-1})$ for $G$-invariant $M$. 
 
  
  
In general we have 
to look at elements  $g$ in $G$ that take the value $0$ for certain \ir characters. Following \cite{DN}  we say that $g$ is {\it non-vanishing} if $\chi(g)\neq 0
$ for every \ir character $\chi$ of $G,$ otherwise we call $g$  {\it vanishing}.
 Vanishing  group elements are of particular interest in the block theory of finite  groups.  We postpone our comments on this matter until Section 4. Here we limit ourselves to the following well known special case. Let $|G|$ denote the order of $G.$ If $p$ is a prime then $|G|_p$ is the $p$-part of $|G|$, that is, $|G|/|G|_p$  is coprime to $p$. The element $g$ in $G$ is {\it $p$-singular} if $p$ divides the order of $g.$

\medskip
\begin{lemma}\label{pr1}  {\rm\cite[Theorem 86.3]{CR}} Let G be a finite group whose order is divisible by the prime  $p$ and  let $\chi$  an \ir character of G. Suppose that  $\chi(1)$, the degree of $\chi,$  is a multiple of   $|G|_p$. 
 Then $\chi(g)=0$ for every $p$-singular element $g\in G$. \end{lemma}
  
In what follow $p$ is a prime. The \ir characters of degree divisible by $|G|_p$ are refered as those of defect 0 (or of $p$-defect 0 
if $p$ is not clear from the context). From this lemma 
and Theorem~\ref{SingCay-2} we obtain the following  general result on singular Cayley graphs.

 \begin{corol}\label{pr2} Suppose that  $p$ divides 
 $|G|$ and that $G$ possesses  an \ir character of   $p$-defect $0$. Then the  Cayley graph ${\rm Cay}(G,H)$ is singular whenever $H$  is a  $G$-invariant connecting set that consists of $p$-singular elements. 
\end{corol}

At the first sight such characters do not appear to be  a common phenomenon. However, this is not so as the following result shows:

\begin{propo}\label{ma13}  {\rm \cite[ Corollary 2]{GO}}   Let $G$ be a non-abelian finite simple group and let $p>3$ be a prime dividing  $|G|.$ Then $G$ has an \ir character of defect $0$. 
 This remains true for $p=2,\,3$ unless $G$ is a sporadic group (with known exceptions) or an alternating group.  \end{propo}
 
It follows by Lemma \ref{pr1}  that  in a non-abelian simple group $G$ any element of order divisible by the prime $p>3$ vanishes at some \ir character of $G.$ Hence Theorem \ref{SingCay-2} applies and combining Proposition \ref{ma13}\, with Theorem \ref{SingCay-2}\, we get

 \begin{theo}\label{tc1}  Let $p>3$ be a prime. Let G be a non-abelian  simple group and $M\subset G\setminus \{1\}$  a G-invariant  subset consisting of p-singular elements. Then the Cayley graph ${\rm Cay}(G,M\cup M^{-1})$ is singular. 
 This remains true for $p=2,3$ unless $G$ is an alternating group or a sporadic simple group.
 \end{theo}
 
  The exceptions in this theorem are genuine. They can be detected easily for alternating groups $A_n$ with $n=7,11, 13$ by inspection of the character tables.
However, for arbitrary $n$ the problem of describing all non-vanishing elements in $A_{n}$ is still open. 
In any case, Theorem \ref{tc1}  yields many examples of singular Cayley graphs.

 We first state   some  elementary results which yield a variety of singular Cayley graphs  when $G=A_{n}$.
  For  $g\in G$ let $(c_1\ld c_k),$ with $c_1\geq c_2\geq \cdots \geq c_k,$ be the cycle lengths of $g,$ in the sense that $g$ has $k$ cycles where the longest cycle is of length $c_{1},$ the second longest of length $c_{2},$ and so on.

  \begin{theo}\label{tc2}  Let $G=A_n$ with  $n>4$ be the alternating group. 
  \\[5pt]
 $(1)$ Let $R_1=\{g\in G: c_k=1,c_{k-1}>1\}$ and let $M\subseteq R_1$ be a $G$-invariant subset. Then $H=M\cup M\up$ is a connecting set  and 
  the Cayley graph  ${\rm Cay}(G,M\cup M^{-1})$ is singular.\\[5pt]
$(2)$ Let $R_2=\{g\in G: c_{k-1}>1$ and $c_i\neq2,4 \}$ for $i=1\ld k$. Let  $M\subseteq R_2$ be G-invariant. Then $H=M\cup M\up$ is a connecting set and 
  the Cayley graph  ${\rm Cay}(G,M\cup M^{-1})$ is singular.
 \end{theo}
 
 Note that $R_1$ consists of all elements of $G=A_n$ fixing exactly one point of the natural  $G$-set, whereas $R_2$
 consists of elements fixing at most one point and having no $2$- and $4$-cycles in their cycle decomposition.

 An element $g\in G$ is called {\it real\,} if $g\up$ is conjugate to $g$. In symmetric groups all elements are real. This is not the case for alternating groups. From Theorem \ref{tc2} we deduce the following 
  
   \begin{theo}\label{nr4a} Let $G=A_n$ with $n\geq 4$ and let 
   $M\subset G$ be any  set of  non-real elements. Then there exists an \ir character of G that vanishes on all   elements of $M$.
Furthermore, $H=M\cup M\up$ is a connecting set and the Cayley graph $Cay(G,H)$ is singular.
\end{theo} 
  
  The proof of Theorem \ref{tc2} is based on the Murnahgan-Nakayama formula for computing the values
  of \ir characters of symmetric groups (see Section 4). In general, the problem of describing, for a given 
  \ir character $\chi$ of a given  group $G,$ the set $\{g\in G: \chi(g)=0\}$ seems to be intractable, even when $G$ is a symmetric or alternating group. However, the block theory of group characters supplies powerful tools for approaching this problem. In particular, we use block theory to prove Theorems \ref{main2} and \ref{pp7} below.

  Recall that for every prime $p$ dividing the order of a  finite group $G$ the set of all \ir characters of $G$ is partitioned into {\it blocks\,} (or {\it $p$-blocks}
  to be accurate)\, and each block $B$ determines a $p$-subgroup of $G,$ defined up to conjugacy in $G.$ This group is the {\it defect group} of $B$. If $\chi\in B$ then $\chi(g)=0$ whenever the $p$-part of $g$ is {\it not} contained in
  a defect group of $B$. To use this fact, it is important to know the blocks with smallest defect groups.  If $p>3$ then this smallest defect group is $\{1\}$ for $A_n,$ see \cite{GO}, and this yields Proposition \ref{ma13} above.  
  (Note that blocks with trivial defect group are called {\it blocks of defect} 0.) For $G=A_{n} $ and $p=2$ the smallest defect groups can be easily determined (see Section  5), whereas for $p=3$ this is still an open problem. This is discussed in 
  \cite[Theorem 2.1]{CCL}, where, for $|G|_3=3^a$ and $n\neq 7$,  
   the order  of a smallest defect group is bounded from above by $ 3^{(a-1)/2}$.
 We improve this bound to  $ 3^{(a-1)/3},$ see Proposition \ref{3au}.

Now we turn to the simplest (in a sense) version of the singularity problem for Cayley graphs: we assume that $H=C\cup C^{-1}$ where $C\neq \{1\}$ is a single conjugacy class in $G$.  Theorem~\ref{tc1}  resolves this version of the problem for Cayley graphs of the shape ${\rm Cay}(G,H)$ with $G$ simple, except when $G=A_{n}$ and when the elements of $C$ have order  $2^\al 3^\beta$ for some $\alpha$ and $\beta.$ 
Theorem \ref{nr4a} reduces the problem to the case where $C=C\up$, that is, where $H$ is a single conjugacy class. 
Below we state some partial results.  One of them   is the following:

\begin{theor}\label{main2}  Let $G=S_n$ or $A_n$ with $n\geq 5$ and $n\neq 7,11,$ and let $\om(G)$ be the set of element orders  of   $G.$   Let $\om_{2,3}(G)$ be the set of all numbers in $\om(G)$ that  are not divisible by   any prime $p>3$. Then $G$ contains a vanishing element of order $m$ for every $1\neq m\in \om_{2,3}(G)$. 
\end{theor}

More can be said about the possible choices of $h$ as an element of order $m.$ By Theorem~\ref{tc1} there are no restrictions unless $m$ is of the shape $2^\al 3^\beta$ for some $\alpha$ and $\beta.$ In that latter case $h$ can be chosen as any element of order $m$ {\it fixing a least number of members} of the natural set $\{1,...,\,n\},$ see Theorem~\ref{mt2}.


  \begin{theor}\label{pp7}  Let $G=A_n$ with $n>4$ and let $g\in G.$ Suppose that $2|g|$ and $3|g|$ are not in $\omega(G)$. 
Then $g$ is vanishing unless $n=7$. 
\end{theor}
  
This statement is new only for   $|g|=2^\al 3^\beta,$ otherwise it follows from  Lemma \ref{pr1} and Proposition \ref{ma13}.
It is not true that all elements satisfying the condition in Theorem \ref{pp7} vanish at the same character of $G.$ But if $\al\beta=0$
then this is the case, see  Corollary  \ref{ma2} and Proposition  \ref{ma3}.

{\sc Notation:} Our notation for finite simple groups agrees with the {\sc Atlas}~\cite{Atl}. In particular, $A_n$ means the alternating group on $n$ letters, and $S_n$ is the symmetric group. The underlying set is often denoted by $\Omega_n$, and it can be idetified with $\{1\ld n\}$.  For a set $M\subset S_n$ the {\it support\,} of $M$ is  $\supp(M):=\{x\in \Omega_n: gx\neq x$ for some $g\in M\}.$ In the other words,   $\supp(M)$ is the complement in $\Omega_n$ of the set of the elements fixed by $M$.

 If $G$ is a group then we write $|G|$ for the order of $G;$ if $p$ is a prime then $|G|_p$
is the $p$-part of $|G|,$ equivalently, the order of a \syl of $G$. For non-zero integers $m,\,n$ we denote  the g.c.d. of $m,n$ by $(m,n).$  If $g\in G$ then $|g|$ is the order of $g$. The identity element of $G$ is denoted by $1$.  For $h\in G$ we write $h^G$ for the conjugacy class of $h$ in $G$. We write $\Irr G$ for the set of all \ir characters of $G.$ If $\chi$ is a character and $M\subset G$ then $\chi(M)=0$ means that $\chi(g)=0$ for all $g\in M$.

\section{ Singularity of Graphs and Cayley Graphs}

Let $\Ga=(V,E)$ be a graph with vertex set $V$ and edge set $E.$ Let $n:=|V|.$ Two distinct vertices $u$ and $v$ are adjacent to each other, denoted $u\sim v,$ if and only if $\{u,\,v\}\in E.$

Let $\Fi$ be a field of characteristic 0. Then we denote by $\Fi V$ the vector space over $\Fi$ with basis $V.$ This is a permutation module for the automorphism group of $\Ga.$ 
The
natural inner product on $\Fi V$ is given by $\lp v,v'\rp=1$ if $v=v'$ and $\lp v,v'\rp=0$ if $v\neq v',$ for $v,\,v'\in V.$ 
The {\it adjacency map} $\alpha\!:\, \Fi V\to \Fi V$ is the linear map given by \begin{equation}\quad\alpha(v):=\sum_{v\sim v'}\,v'\end{equation} for $v\in V.$ 
Since $v\sim v'$ for $v,\,v'\in V$ if and only if $v'\sim v$ we have $\lp \alpha(v),v'\rp=\lp v,\alpha(v')\rp.$ Therefore $\alpha$ is symmetric with respect to  this inner product.  The matrix of $\alpha$ with respect to the basis $V$ is  the {\it adjacency matrix} $A=A(\Ga)$ of $\Ga.$ Since $A$ is symmetric $A$  is diagonalizable when $\Fi=\re$ and so
all eigenvalues are real. The {\it spectrum} of $\Ga$ are the eigenvalues $\gl_{1},\,\gl_{2},\,...,\,\gl_{n}$ of $A.$

The graph $\Ga$ is {\it singular} if $A(\Ga)$ is a singular matrix. In other words, $\Ga$ is {singular} if and only if $0$ is an eigenvalue of $\Gamma.$ Singularity is connected to another special eigenvalue of  graphs. The {\it complement} $\bar{\Ga}$ of $\Gamma$ is the graph on the same vertex set with two distinct vertices connected in $\bar{\Ga}$ if and only if they are not connected in $\Ga.$ For regular graphs (each vertex has the same number of neighbours) one can show easily that $\Ga$ is singular if and only if $-1$ is an eigenvalue of $\bar{\Ga}.$ The {\it nullity ${\rm null}(\Ga)$} of $\Ga$ is the nullity of $A(\Ga).$ Singular graphs are the graphs with ${\rm null}(\Ga)>0.$ For instance, if $\Ga$ is a bipartite graph with parts $V=V_{1}\,\dot{\cup}\,V_{2},$ then  ${\rm null}(\Ga)\geq \big||V_{1}|-|V_{2}|\big|.$ This observation  provides a wealth of examples of singular graphs.

\medskip
The singularity of graphs plays a significant role in several parts of mathematics and applications. It would be impossible to review the vast literature on graph eigenvalues in this paper. In representation theory and finite incidence geometry the containment of one permutation character in another often is easiest to establish by showing that a certain graph is non-singular or that its nullity is bounded in a particular way. A famous example is the theorem of Livingstone and Wagner about the representations of a permutation group $G$ on the $k$- and $(k+1)$-subsets of the set on which $G$ acts.   

\medskip
We mention also the significance of graph singularity in systems analysis, physics and chemistry, see for instance the survey article \cite{GutBo}. Essentially, when modelling a discrete mechanical system (Hamiltonians) 
it is often necessary to work out a  linear approximation of an operator where the constituents of the system and the relationships between them are represented by a finite  graph. Many characteristics and observables of the system -- its energy for instance -- then typically involve the spectrum of this graph. This is one of the  principles that underpins spectroscopy and  H\"{u}ckel Theory in chemistry~\cite{AStreit}. In such applications the singularity of a molecular graph of a feasible compound typically indicates that the compound is highly reactive,   unstable, or nonexistent, see \cite{GutBo}.  


\bigskip
In this paper we concentrate on the singularity of graphs whose automorphism group is transitive on vertices. This includes in particular Cayley graphs for which we now give the basic definitions. 

\medskip 
Let  $G$ denote a finite group with identity element $1.$ Then the subset $H$ of $G$ is  a {\it connecting set} provided the following holds: \\[10pt]   (i) $1$ does not belong to  $H,$\\[5pt] (ii)
$H^{-1}:=\{h^{-1}\,|\,h\in H\}=H$\, and\\[5pt] (iii) $H$ generates $G.$ 

Suppose now that $H$ is a connecting set. Then define the graph $\Ga=(V,E)$ with vertex set  $V=G$ by calling two vertices $u$ and $v\in G$  {\it adjacent}, denoted $u\sim v,$ if  there is some $h$ in $H$ with $hu=v.$ The first condition above is equivalent to saying that $\Ga$ has no loops. The second conditions holds if and only if all  edges are undirected, that is $u\sim v$ if and only if $v\sim u.$ The last condition is equivalent to saying that $\Ga$ is connected. This graph is the {\it Cayley graph $\,\Ga = {\rm Cay}(G,H)\,$ on $G$ with connecting set $H.$} Its  {adjacency map} $\alpha\!:\, \Fi V\to \Fi V$  has the form \begin{equation}\quad\alpha(v):=\sum_{h\in H}\,h^{-1}v\end{equation}  for all vertices $v$ in $V=G.$ 
Since $H=H^{-1}$ the set of all neighbours of $v\in \Ga$ is the set $Hv.$ In particular, $\Ga$ is regular of degree $|H|.$ Similarly, $HHv$ is the set of all vertices of distance $\leq 2$ from $v,$ and so on. The radius $r(\Ga)$ of $\Ga,$ as a graph invariant, is useful for studying generating sets in a group. Evidently, $r=r(\Ga)$ is the least number $r>0$ such that $H^{r}:=\{h_{1}h_{2}\cdots h_{r}\,|\, h_{i}\in H\}$ is equal to $G.$ This invariant is a subject of intensive study by group theorists.

Let $g\in G.$ Then the right multiplication $x\mapsto xg$ for $x\in G$ is an automorphism of $\Ga$ as is easy to see. Therefore the  
right-regular representation of $G$ on itself provides an injective homomorphism $G\to {\rm Aut}({\rm Cay}(G,H))$ for any connecting set $H.$ Cayley graphs are characterized by this property:

\medskip
\begin{theor}\label{Sabi}(Sabidussi
) The graph $\Ga$ is isomorphic to 
 a Cayley graph 
 if and only if  ${\rm Aut}(\Ga)$ contains a subgroup that acts regularly on the vertices of $\Ga.$
\end{theor} 

Next consider  the left  multiplication $x\mapsto g^{-1}x$ for $g$ and $x\in G.$ By contrast,  this does {\it not\,} yield an  automorphism of $\Ga$  in general.  It is  easy to show that  $x\mapsto g^{-1}x$ is an automorphism of ${\rm Cay}(G,H)$ if and only if $gH=Hg.$ This is relevant for this paper as we are  dealing with connecting sets that are unions of conjugacy classes. 

If $\Gamma={\rm Gay}(G,H)$ then the space $\Fi G$ is the underlying space of the group algebra of $G$ over $\Fi$.  Then  
    $ \gaa=\sum_{h\in H}\,\rho(h)$ is an element of the group algebra, which is in the center of it whenever $H$ is a union of conjugacy classes.
 
\medskip  
{\it Proof of Theorem} \ref{SingCay-2}: Let $\rho$ be the left regular representation, $\rho(g)(v)=g^{-1}(v)$ for $g\in G$ and $v\in V=G.$ The adjacency map  (2)  above then becomes  $$\quad \gaa(v)=\lp\sum_{h\in H}\,\rho(h)\rp(v)\,.$$
Note that $\rho(g)\gaa=\gaa\rho(g)$ for $g\in G$ as  $g(\sum_{h\in H} \rho(h))g^{-1}=\sum_{h\in H} \rho(h)$.  

Let $\gl_{1},\,\gl_{2},...,\,\gl_{t}$ be the distinct eigenvalues of $\gaa$ and let $E_{1},\,E_{2},\,...,\,E_{t}$ be the distinct eigenspaces of $\gaa.$ So $\Fi V=E_{1}\oplus E_{2}\oplus\,...\oplus E_{t}$ and each $E_i$ $(i=1\ld t)$ is invariant under $\rho(g)$ for every $g\in G$. Let $\rho_{i}$ and $\gaa_{i}$ denote the restriction of $\rho$ and $\gaa$ to $E_{i},$ respectively. Thus $$\gaa_{i}(x)=\lp\sum_{h\in H}\,\rho_{i}(h)\rp(x)=\gl_{i} x\,$$ for $x\in E_{i}.$ 
Now, if ${\rm Cay}(G,H)$ is singular, say $\gl_{1}=0,$ then every irreducible representation $\rho_{1,i}$ appearing in $\rho_{1}$ satisfies $\sum_{h\in H}\,\chi_{1,i}(h)=0,$ where $\chi_{1,i}$ is the character of $\rho_{1,i}.$ Conversely, if $\chi_{j,i}$ is an irreducible character with $\sum_{h\in H}\,\chi_{j,i}(h)=0$ then $\rho_{j,i}$ appears in $\rho$ and so there is some $E_{j}$ on which $\gl_{j}=0.$ \dne
  
Finally we consider a connected graph $\Ga=(V,E)$ which admits a vertex transitive group $G$ of automorphisms. In this case we construct an associated Cayley graph $\Ga^{*}:={\rm Cay}(G,H)$ as follow.  Fix a vertex $v\in V$ and let $C$ be its stabilizer in $G,$ with $c:=|C|$. In view of Sabidussi's theorem we may assume that $c>1.$  Next let $H:=\{\,h\in G\,:\, v\sim v^{h}\,\}.$ Clearly $1\not\in H$ and $H=H^{-1}.$ Also, $H$ generates $G,$ this follows from the transitivity of $G$ on vertices and the connectedness of $\Ga.$ Therefore we have a Cayley graph $\Ga^{*}:={\rm Cay}(G,H)$  associated to $\Ga.$ It is {\it imprimitive\,} in the sense that $a'\in Ca$ is adjacent to $b'\in Cb$ in $\Ga^{*}$ if and only if $v^{a}$ is adjacent to $v^{b}$ in $\Ga$, for all $a,a',b,b'\in G.$ Hence the adjacency matrix of $\Ga^{*}$ if of the  form $A(\Ga)\otimes J$ where $J$ is the $c\times c$-matrix with all entries equal to $1.$ This implies the following

\begin{theor}\label{TransGraph}{\rm \cite{Lo}} Let $\Ga=(V,E)$ be a connected graph with a vertex transitive group $G$ of automorphisms and let $\Ga^{*}$ be the associated Cayley graph. Suppose that the eigenvalues of  $\Ga$ are $\gl_{1},...,\gl_{n},$  with $n=|V|.$ Let $c$ be the order of  the stabilizer in $G$ of a vertex of $V.$ Then the eigenvalues of  \,$\Ga^{*} $ are  $c\gl_{1},..,\,c\gl_{n}$ together with $0,...,\,0,$ of multiplicity  $n(c-1).$  In particular, $\Ga$ is singular if and only if ${\rm null}(\Ga^{*})>n(c-1).$
\end{theor}


{\sc Comment:\,} We see that the singularity problem for vertex transitive graphs can be reduced - in principle at least - to the nullity problem for Cayley graphs.  
The theorem can also be used to construct singular graphs: any graph with a vertex transitive but not vertex regular group of automorphisms yields a singular Cayley graph with the same group of automorphisms. 

\def\cc{character }

\section{Elementary observations on zeros of alternating groups characters }

The comments in Section 1  suggest to pay particular attention to the alternating groups. In fact, the reasonings in this paper are mostly concerned with these groups.  In this section we collect a number of well known facts about characters of alternating groups and prove some results on the zeros of some of their \ir characters. 

 We
  first recall certain notions of the \rep theory of $S_n$.  It is well known that the \ir characters of $S_n$ are in bijection with the Young diagrams, and also with the partitions of $n$. So we write $\phi_Y$ for the \irr  or the  \ir character of $S_n$ corresponding to the Young diagram $Y$. 
 For the  Young diagram $Y$ we write $|Y|$ for the number of boxes in it.  A subdiagram of   $Y$ is a Young diagram of $S_m$ for $m<n$ which is contained in $Y$ as a subset with the same top left hand corner.  A box in  $Y$ is called extremal if there is no box either below or to the right of it. The set of all extremal boxes form the {\it rim} of $Y.$

The notion of a hook in a Young diagram $Y$ is common knowledge, see \cite[page 55]{JK}. The number of boxes in a hook
is called the {\it length} of it. A  hook of length $m$ is called an $m$-{\it hook}. The {\it leg} of a hook is the set of all  boxes below the first row and ends in its {\it foot}. The number of the boxes in the leg is the {\it leg length.} The {\it arm} of the hook is its horizontal part, it ends in the {\it hand} of the hook, the right furthest box in the arm.  Both foot and hand of the hook belong to the rim of $Y.$
 
 To every $m$-hook $\Gamma$ of $Y$ there corresponds the set $\nu(\Gamma)$ 
 of $m$ contiguous boxes lying on the rim of $Y$ which link the foot  to the  hand of $\Gamma.$ These boxes are the {\it hook rim} of $Y,$ see \cite[pages  56 and 75]{JK} for details. Removing $\nu(\Gamma)$  from $Y$ yields a Young subdiagram $Y\setminus \nu(\Gamma)$ with 
$|Y|-m$ boxes. 
Conversely, if $\nu$ is a sequence of $m$ contiguous boxes on the rim of $Y$ then the ends of $\nu$ are the foot and hand of a unique hook $\Gamma$ with  $\nu=\nu(\Gamma).$ We call $\nu$ an $m$-rim of $Y$ if $\nu=\nu(\Gamma)$ for some hook $\Gamma$ of $Y$, and refer to the leg length of $\Gamma$ as the {\it leg length} of $\nu. $  
 

\medskip

Below we need the  Murnahgan-Nakayama formula \cite[2.4.7]{JK}. It expresses the character value of an \ir character $\chi_Y$ in combinatorial terms. Let $g\in S_n$ and $g=ab$ where $a$ is an $m$-cycle and where $b\in S_{n-m}\subset S_n$ is the permutation induced by $g$ on the points fixed by $a.$ The Murnahgan-Nakayama rule is the induction formula 
\begin{equation}\label{eq MN}
\chi(g)=\sum (-1)^i\chi_{(Y\setminus \nu)}(b)\end{equation} where the sum runs over all $m$-rims $\nu$ of $Y$ and where  $i$ is  the leg size of $\nu$. (If no  $m$-rim  exists then we have $\chi(g)=0$ by convention.)




\medskip
As an illustration, we state the following



  \bl{n44}  Let $G= A_n$ or $S_n$ with $n\geq 7$  and let $M\subset  G$ be the subset of all elements whose cycle decomposition has a cycle of length greater than $2\sqrt{n}+2$. Then 
  $\chi(M)=0$ for some \ir character   $ \chi$ of G.
  \el
  
\pf Let $m$ be the minimal number $i$ such that $i^2>n$, so $m>\sqrt{n}$. Let $g\in M$ and let $c(g)$ be maximal length 
  of a cycle in the cycle decomposition of $g$. Then $c(g)>2  \sqrt{n}+2$.  If $m(m-1)>n$ then set $Y=[m\ld m,n-m(m-1)]$ ($m$ rows),
  if $m(m-1)=n$ then set $Y=[m\ld m]$ ($m-1$ rows), if $m(m-1)< n$ then set $Y=[m-1\ld m-1,n-(m-1)^2]$ ($m$ rows).
 (So $Y$ is a nearly square diagram.)  In all cases the hook lengths of $Y$ does not exceed $2m-1$.  By the Murnahgan-Nakayama formula  (\ref{eq MN}) 
we have $\chi_Y(g)=0$ whenever $c(g)\geq 2m.$ As  $c(g)>2  \sqrt{n}+2\geq 2(m-1)+2$, 
  the result follows. \dne 

  
 

   \begin{corol}\label{ma2}  Let $ G=A_n$ with $n>4$  and let $e$ be the maximum order of a  $2$-element   of $G.$   Let M be the set of elements of order $e.$ 
 Then $\chi(M)=0$ for some \ir character $\chi$ of $G$. 
 \end{corol}
 
\pf Let $g$ be an element of $M.$ Observe that $| g|\geq n/2$. (In $S_n$ the inequality is strict.) If $\frac{n}{2}>2\sqrt{n}+2$ then the result follows from    Lemma \ref{n44}. Let  $\frac{n}{2}\leq2\sqrt{n}+2$, equivalently,  $n^2-24n+16\leq 0$. This implies $n\leq 23$. Let $c$ be the maximal cycle length in the cycle decomposition of $g$. Then   $c=16$ for $18\leq n\leq 23$
and  $16>2\sqrt{n}+2$  so Lemma~\ref{n44} applies. If $  n= 6,8,10,12,15,17$ then $G$ has a character of defect 0 which vanishes at all 2-elements of $G,$ by  Lemma~\ref{go1}. 
 If $n=16, 14,13,11$ then $c=8.$ Let $Y$ be as in the proof of Lemma~\ref{n44}. Then the maximal hook length of $Y$ is at most $7$ and 
 so the result follows as above by the Murnahgan-Nakayama formula  (\ref{eq MN}). If $n=7, 9$ then  the lemma follows by inspection of the character table of $G$. \dne


  \bl{nr5}  Let $g\in A_n$ with $n\geq 7$  and let   $(c_1\ld c_k)$ with $c_1\geq c_2\geq \cdots\geq c_k$ be the cycle lengths of $g.$   Let $\chi$ be the character labeled by the Young \di  $Y=[n-4,3,1]$. Suppose that $c_{k-1}\neq 1$ and $c_i\neq 2,4$ for $i=1\ld k$. Then 
  $\chi(g)=0$. 
   \el


\pf  For the reader's convenience we visualize the shape of \,$Y=[15,3,1]$\, for $n=19$:
   \bigskip

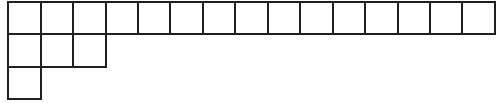
\centerline{ $ \yng(15,3,1)$   }

   We view $\chi$ as a character of $S_n$.  Let $g=g_1b$ where $g_1$ is a cycle of size $c_1$ and the cycle lengths of $b$ are $(c_2\ld   c_k)$.
   By the Murnahgan-Nakayama rule (\ref{eq MN})
    $\chi(g)=\sum (-1)^i\chi_{(Y\setminus \nu)}(b)$ where $\nu$ runs over the
    $c_1$-rims of $Y$ and  $i$ is the leg length of $\nu$. (If no  $c_1$-rim exists then $\chi(g)=0$.)
Set $r=c_1$. 
 It is clear from the \di shape that an  $r$-rim is either a part of the first row
(and then $r\leq n-7$) or $Y_1=Y\setminus \nu$ is one of the diagrams $[2,2,1]$, $[2,1^2]$ or $[2]$. In each case there is 
exactly one way to delete $\nu,$ so the sum has at most one term. 

Suppose first that $n-r<7$. There is no way to remove  an $r$-rim to obtain $Y_1$ of size 6. So $\chi(g)=0$ if $n-r=6$. 
Let $n-r\leq 5$; then
$c_2+\cdots +c_k=n-r\leq 5$,  and hence  $c_2=|b|=5,3$ or 1. The options $b=1$ and $|b|=3$ for $Y_1=[2,2,1]$ are ruled out as $c_2\neq 1$. 
Let $\chi_1$ be the character of $S_{n-r}$ corresponding to $Y_1$. Then $\chi_1(1)=5,3$ for $Y_1= [2,2,1]$, $[2,1^2]$, respectively. Therefore, $\chi_1$ is of $c_2$-defect 0, and hence $\chi_1(b)=0$ (one can also check this in the character tables of $S_5,S_4$). This implies $\chi(g)=0$ in these cases.

Let $n-r\geq 7$. Then $Y_1=[n-r, 3,1]$ and the leg length of the  $r$-rim $\nu$ is 0. So $\chi(g)=\chi_1(b)$. So we can use induction on $k$. The case with $k=1$ follows from the above as then $n-r=0<7$. If $k>1$ then  the cycle lengths  of $b$
are $(c_2\ld c_k)$, so the result follows by the induction assumption. \dne

 \begin{propo}\label{k9}  Suppose that $G$ is a doubly transitive permutation group on $\gO.$ Then every element $g\in G$ fixing exactly one point of $\gO$ is vanishing. In particular, if $G=A_{n}$ with $n>3$ and if $g$ fixes exactly one point then $g$ vanishes at the \ir character of degree $n-1$.
 \end{propo}
 
 \pf The permutation character is of the shape $\pi=1+\chi$ with $\chi$ irreducible. By assumption therefore $\chi(g)=\pi(g)-1=0.$ \dne

{\it Proof of Theorem } \ref{tc2}. This follows by combining Theorem \ref{SingCay-2} with Lemma \ref{nr5} and Proposition \ref{k9}.

\bl{re1}
 An element  $g\in A_n$ is a real element  if and only if  one of the following conditions holds:\\[5pt]
$(1)$ The cycle decomposition of $g$  has a cycle of even length.\\[5pt]
$(2)$ The cycle decomposition of $g$ has two cycles of equal odd length. Note that the fixed points are counted as cycles of length $1$, so this includes any permutation that has two or more fixed points.\\[5pt]
$(3)$ All the cycles of $g$ have distinct odd lengths  $c_1\ld c_k$ and $\sum_{i=1}^k(c_i-1)/2$  is even. In other words, the number of $c_i$'s that are congruent $3$ to  modulo $4$  is even.\el

\pf Clearly $g$ is conjugate to $g^{-1}$ in $S_{n}.$ Note that the $S_{n}$-conjugacy class of $g$ is an $A_{n}$-conjugacy class if and only if there is an odd permutation that centralizes $g,$ these are the conditions 1) and 2).  In the remaining case, a cycle of odd length  $2\ell+1$ is inverted by an element of sign $(-1)^{\ell},$ and this gives the condition 3). Here the $S_{n}$-conjugacy classes  of $g$ split in $A_n$, but  $g$ is conjugate to $g^{-1}$ in $A_{n}.$\dne

\medskip
{\it Proof of Theorem} \ref{nr4a}. If $n\geq 7$ then any non-real element satisfies the assumption of Lemma \ref{nr5}, whence the result. If $n\leq 6$ then either $n=6$, $|g|=5$ or $n=4$, $|g|=3$. In these cases the result follows from Lemma \ref{k9}. 
Finally, the claim that $H$ is connected follows immediately if $n>4$ as $A_n$ is simple. If $n=4$ then $H$ consists of all elements of order 3, so the claim follows by inspection of normal subgroups of $A_4$.\dne

  \begin{propo}\label{nr7} Let G be a simple group and suppose that $g\in G$ is non-real. Then $g$ is vanishing.
 \end{propo}
  
\pf  By  Theorems \ref{tc1} and \ref{nr4a},  we are left to check the sporadic groups. It is observed in \cite[p. 414]{INW} that 
  $M_{22}$, $M_{24}$ are the only sporadic groups having non-identity non-vanishing elements, and these are of order 2, see the {\sc Atlas}~\cite{Atl}, and hence are real. \dne

 \bl{sd9}  Let $G=S_n$ with $n>2$ and let $M\subset  (S_n\setminus A_n)$ be a subset.
 Then $\chi(M)=0$ for some \ir character of G. Furthermore, if $M$ generates G then the Cayley graph ${\rm Cay}(G,M\cup M\up)$ is singular.
 \el
 
\pf  Note that there exists a symmetric Young diagram of size $n$. Let $\phi$ be 
the \irr of $G$ labeled by this diagram. It is well known that $\phi$ is reducible under restriction to $A_n$ and that the irreducible constituents are non-equivalent. Let $h\in M$. By Clifford's theorem, $h$ permutes these constituents. So if $V$ is the underlying space of $\phi$ then $V=V_1\oplus V_2$, and $hV_1=V_2$, $hV_2=V_1$. It is clear from this that the trace of $h$ equals $0$.
So the character of $\phi$ vanishes on $M$. If  $M$ generates $ G$  then the set $H=M\cup M\up$ is connecting, so  the result follows by Theorem \ref{SingCay-2}.  \dne

 \section{Blocks in symmetric groups and vanishing elements}
 
 In this section we expose some part of representation  theory of symmetric groups that is needed for the remainder of the paper. Recall the notation from the end of Section 1.
 Let $G$ be a finite group. For every prime $p$ dividing $|G|$ the \ir characters of $G$ partition into $p$-blocks. To every $p$-block there corresponds a conjugacy class of $p$-subgroups of $G$ and each of them is  called a {\it defect group}\, of the block. If $p^d$ is the order of  a defect group then $d$ is called the {\it defect of the block}.  In particular, blocks of defect 0 are those whose   defect groups consist  of one element.   In addition, $G$  has a {\it $p$-block of defect $0$} if and only if there is an irreducible character of $G$ which has $p$-defect $0.$  See for instance Navarro~\cite{na1} or Curtis and Reiner~\cite{CR}, Chapter VII, for general theory of blocks.   We  use the following well known fact:

\begin{lemma}\label{jk1}  {\rm \cite[Corollary 15.49]{Is1}}  Let G be a finite group and let $g\in G$ be a p-singular element. Let $g=g_ph=hg_p$, where $g_p,h\in\lan g \ran$, $g_p$ is a p-element and $|h|$ is coprime to  $p$.  Let $\chi$ be an \ir character of G.
 Suppose $g_p$ is not contained in any defect group 
 for the $p$-block containing $\chi$. Then $\chi(g) =0$, in particular $g$ is a vanishing element of $G.$
  \end{lemma}
 
To use  this lemma, one needs to know the defect groups of the $p$-blocks of $S_n$ (for $p=2$ and 3). These are described in
 \cite[Theorems 6.2.39 and 6.2.45]{JK} for any prime $p$. However, first we discuss a special case of blocks and characters of defect 0.

\subsection{Characters of defect $0$}

For  non-abelian simple groups  there is the following criterion for the existence of $p$-blocks of defect $0.$

\begin{lemma}\label{go1}  {\rm \cite[ Corollary 2]{GO}} Every non-abelian finite simple group $G$ has a $p$-block of defect $0$ for every prime
$p,$ except in the following special cases:\\[5pt]
{\rm (1):}\,\, $G$ has no $2$-block of defect $0$ \ii  G  is isomorphic to $M_{12}$; $M_{22}$; $M_{24}$; $J_2$; $HS$; ${\rm Suz};$ ${\rm Ru};$ $C_1$; $C_3$; {\rm BM}, or $A_n$ where $n \neq 2m^2+m$ and  $n\neq 2m^2 +m + 2$ for any integer $m$ (not necessarily positive).\\[5pt]
{\rm (2):}\,\, $G$ has no $3$-block of defect $0$ \ii  G is isomorphic to ${\rm Suz},$ $C_3$, or $A_n$ with $3n+1 =m^2r$, where r is square-free and divisible by some prime $q \equiv 2 \pmod 3$.
 \end{lemma}

{\sc Remarks:} \,1. For $G=A_n$ the condition (1) can be expressed in an alternative way. If $m\geq 0$ then set $k=2m+1$. Then $2m^2+m=\frac{k(k-1)}{2}$. 
 If $m<0$ then set $k=-2m$. Then $2m^2+m=\frac{k(k-1)}{2}$. So the set $\{2m^2+m:m\neq 0\}$ coincides with the set 
  $\{\frac{k(k-1)}{2}:k>0\}$. So $A_n$  has a 2-block of defect 0 \ii $n=\frac{k(k-1)}{2}$ or $2+\frac{k(k-1)}{2}$ for $k>0$.\\[5pt]
 2.\, For $p=3$ the condition (2) can be expressed in an alternative way. Specifically,  $A_n$ has a 3-block of defect 0 \ii $n$ is of the form $n=
 3(x_1^2+x^2_2 +x_1x_2)+x_1+2x_2$, where $x_1,x_2$ are integers, not necessary positive, see \cite{K} or \cite[p.333]{GO}.


 \begin{corol}\label{i51}  If $n\equiv 3\pmod 4 $ then $A_n$ has no $3$-block of defect $0.$
 \end{corol}
  
\pf  Write $n=4l+3$. Then $3n+1=2(6l+5)$, and the number $6l+5$ is odd. If we write $3n+1=2(6l+5)=m^2r$ then $m$ is odd and $r=2r'$, where $r'$ is odd and square free. \itf  $2r$ satisfies the condition in Lemma \ref{go1}(2), whence the result. \dne
 
 Applying Lemma \ref{go1} to a specific number $n,$ and decomposing $3n+1$ as a product of primes,  we obtain the following list of  $n<60$ with $n\not\equiv 3\pmod 4$ for which $A_n$ does not have a $3$-block of defect 0:  $n=13, 18,28,29,38,45,46,48,53,59$.

 \medskip
 {\sc Remark:} \,It is not true that  if  $G$ does not have an \ir character of $3$-defect $0$ then  there exists a $3$-singular element  $g\in G$ such that  $\chi(g)\neq 0$  for every \ir character   $\chi$ of $G$, see the character table of $Suz$.

\begin{corol}\label{go2}  Let $G$ be a perfect 
group, $h\in G$ and $C=h^G$ be the conjugacy class of h.   Suppose that $C$ generates $G$.  
   Then the Cayley graph ${\rm Cay}(G,C\cup C\up )$ is singular except  possibly every simple quotient of G is isomorphic to $A_n$,
 $M_{22}$ or $M_{24}$  and  the order of $h$ is of the form $2^a3^b$ for some integers $a,b.$  
\end{corol}

Recall that  every finite group $G$ has a unique maximal normal nilpotent subgroup $F(G),$ called the Fitting subgroup of $G$.  Lemma~\ref{go1} can be extended to non-simple groups as follows. 

 \begin{propo}\label{dn1} {\rm \cite[  Theorem A]{DN}} Let G be a finite group, and let $h \in G$  be  of order coprime to $ 6$. Then
either $h$ belongs to the Fitting subgroup  $F(G)$ of $G$ or $h$ is vanishing. 
 \end{propo}

For Cayley graphs we therefore have the following general result:

\begin{theo}\label{dh2}  Let G be a non-cyclic  finite  group. Let $1\neq h\in G$ has order  coprime to $6$ and put $H:=h^G\cup (h\up)^G.$   Suppose that $H$ generates $G.$ Then the Cayley graph ${\rm Cay}(G,H)$ is singular. 
\end{theo}


\pf  Suppose first that $G\neq F(G)$. Then $h\notin F(G)$. By Proposition~\ref{dn1}, $\chi(h)=0$ for some \ir character $\chi$ of $G.$ Recall that $\chi(h)=0$ implies $\chi(h\up)=0$. So the result follows by Theorem~\ref{SingCay-2}. Next, suppose that  $G=F(G).$  Then  all   non-vanishing elements  of $G$ belongs to $Z(G)$ \cite[Theorem B]{INW}. 
 Then $G=\lan H\ran$ implies that $G$ is cyclic. \dne

 We mention the following recent result on zeros of characters in $S_n$. This can be used to construct Cayley graphs of relatively high nullity.
 
 \begin{lemma}\label{ma1}  {\rm \cite[Theorem 4.1]{MNO}} Let $p$ be a prime,  let $n\geq p$ be a natural number and $n=a_0+a_1p+\cdots +a_kp^k$ its  p-adic expansion. Let  $h$ be  a  $p$-element of  $S_n$ whose cycle structure  is $1^{a_0}p^{a_1}(p^2)^{a_2}... (p^k)^{a_k}$.  Let $\chi$ be an \ir character  of $S_n$ such that $p$ divides $\chi(1)$. Then $\chi(h)=0$. This remains true for $A_n$ provided $h\in A_n$. $($See Remark following {\rm \cite[Th 4.2]{MNO})}.\end{lemma}

 \subsection{Blocks of Symmetric groups}
 Let $p$ be a prime. Every diagram that does not contain a  $p$-hook is called a $p$-{\it core}. For instance, 2-cores are  the diagrams of triangle shape  $[k,k-1\ld 1];$ in particular, $S_n$ has no $2$-core of size $n$ unless  $n=1+...+k=k(k+1)/2$ for some integer $k>0$.    Every diagram $Y$ contains a unique $p$-core subdiagram  $\tilde Y$ which is  maximal subject to condition $|Y|\equiv |\tilde Y|\pmod p$. The key result of block theory of symmetric groups  states that two \ir characters are in the same block if their Young diagram yield  the same $\tilde Y$  \cite[6.1.21]{JK}.  There is a simple algorithm to obtain $\tilde Y$ as follows.

If $Y$ has no $p$-hook then $Y=\tilde Y$. Otherwise remove arbitrary $p$-rim to obtain a subdiagram $Y_1$. 
If $Y_1$ has a $p$-hook,  remove some $p$-rim from $Y_1$ to obtain a subdiagram $Y_2$ and so on.     The process stops \ii one gets a subdiagram $\tilde Y$  which is a $p$-core.
By \cite[Theorems 2.7.16]{JK}, this final subdiagram $\tilde Y$ is unique (independently from the $p$-hooks choice), and called the $p$-{\it core of} $Y$. Thus,
$|\tilde Y|=|Y|-pb$ for some uniquely determined  integer $b\geq 0$, and this $ b $ is called the $p$-{\it weight of} $Y$ (see \cite[p. 80]{JK}).     Note that the $p$-weight of a diagram is $0$ \ii the diagram is a $p$-core.

The following well known fact follows easily from the dimension formula for \ir characters in terms of hooks:

\bl{d01}  Let $\chi$ be an \ir character of $S_n$ labeled by a Young \di Y. Then $\chi$ is of p-defect $0$ \ii 
Y is a p-core.
\el

 
 
 
 
 
 

By \cite[6.1.35 and 6.1.42]{JK}, there is a bijection between $p$-blocks of $S_n$ and the $p$-cores $C$ such that $|C|\leq |n|$ and $n-|C|\equiv 0\pmod p$.


\begin{theo}\label{t42} {\rm \cite[Theorems 6.2.39]{JK}} Let $\chi$ be an \ir character of $S_n$ labeled by the Young \di $Y,$ and let $B$ be the $p$-block to which $\chi$ belongs. Let $b$ be the p-weight of $Y$. Then a \syl of $S_{pb}$ is a defect group of $B$. 
\end{theo}


Recall that the defect groups of a block are unique up to conjugacy. 
Here  the group $S_{pb}$ is a natural subgroup of $S_n$ in the sense that  this permutes $pb$ elements of $\Omega_n$, and fixes the remaining elements. (Note that if $b=0$,  then $S_{pb}$ is meant to be the identity group, and if $n=pb$ then $S_{pb}=S_n.$) Moreover,  the  character  of $S_{n-pb}$ corresponding to
 $\tilde Y$ is  of defect 0.  

  \begin{corol}\label{uc3}Let D be a defect group of a p-block of $S_n$ and $g\in S_n$ a p-element. Then g is conjugate to an element of D \ii $|\supp(g)|\leq pb=|\supp(D)|$.  
  \end{corol}
 
 \begin{corol}\label{c43} $(1)$ Let B be a p-block of $S_n$ with defect group D and p-core C. Then 
 D fixes exactly $|C|$ elements of $\Omega_n$.
  \end{corol}



It is easy to construct \ir characters with given $p$-core $C$ (provided $n-|C|$ is a multiple of $p$):

 \begin{corol}\label{cc1}   Let $C=[l_1\ld l_k]$ be a p-core and $|C|=c$.  Let $Y=[l_1+bp,l_2\ld l_k] $ be a diagram of $S_{c+bp}$,  $\chi$ the character labeled by Y, and let $ D$ be a defect group of the p-block  the character  $\chi$ belongs to. Then $C=\tilde Y$, and $D $ is a \syl of $S_{pb}$, which is the stabilizer of $c$ elements of $\Omega$ in $S_n$. In particular,  $D$ stabilizes $c$  elements of $\Omega_n$. %
 \end{corol}
  
   For our purpose we are  interested in the defect groups rather than in blocks themselves.  Moreover, we can fix a defect group in every conjugacy class of defect groups 
in such a way that these defect groups form a chain with respect of inclusion. In fact, if the defect groups $D ,D '$ are Sylow $p$-subgroups in $S_{pb},S_{pb'}$, resp., and $b<b',$ then we can assume $D \subset D '$. (For this, one can order the elements of $\Omega$ and choose $S_{pb}$ to be the subgroup
fixing elementwise the last $n-pb$ elements of $\Omega$.) Therefore, with this ordering of defect groups  it is meaningful to speak of the minimal defect group of $S_n$, that is, the one with least possible $b$. Recall that the defect groups of a block are conjugate, and if $D$ is one of them then the defect of a block is the number $d$ such that $|D|=p^d$. So a minimal defect group is a defect group of a block of minimal defect.  Note that the maximal defect group is always a \syl of $S_n$. By Lemma \ref{go1}, if  $p>3$ then the minimal defect group is  trivial (that is, the group of one element).  (Formally, the lemma is stated for  $A_n$ but it remains  true for $S_n$.)

\medskip
If $p>2$ then the defect groups of $A_n$ are exactly the same as those of $S_n$; if $p=2$
then the defect groups of $A_n$ are of shape $D\cap A_n$ for a defect group $D$ of $ S_n$.
  Moreover, if $\chi$ is an \ir character of $S_n$ reducible as that of $A_n$ then
\ir constituents belong to blocks whose defect groups are conjugate in $S_n$, and hence have the same support. (This follows from  \cite[9.26, 9.2]{na1}.)

\bl{ns7}  Let B be a $2$-block (resp., $3$-block)  of $S_n$ of non-zero defect. Then B contains an \ir character that remains 
irreducible under $A_n$.  
\el

\pf It is well known that the characters labeled by non-symmetric diagrams are \ir under $A_n$. So we show that $B$ contains a character 
whose Young diagram is not symmetric. 

  Let  
 $Y=[l_1\ld l_k]$ be the 2-core (resp. $3$-core) diagram determined by  $B$. As $B$ is not of defect 0, $n\neq |Y|$, so $n-|Y|=2b $        (resp. $3b$)  for some integer $b>0$.  If $l_1\geq k$ then the diagram   
 is not symmetric and   the character labeled by $Y_1$  belongs to the block $B$ (see Lemma \ref{cc1}). If $l_1<k$  
 (so $p=3)$ then take for $Y_1$ the diagram  obtained from $Y$ by adding $3b$ boxes to the 1st column, and conclude similarly.   \dne

\bl{23a}  Let $ G=S_n$ or $A_n$, $g\in G$,  and let $D_2$ and $D_3$ be defect groups of a $2$-block, or a $3$-block, respectively. Suppose that $|g|=2^\al 3^\beta$ and $|\supp(g)|>|\supp(D_2)|+|\supp(D_3)|$. Then $\chi(g)=0$ for some \ir character $\chi$ of G. 
\el

\pf  Let $G=S_n$ and  
let $g=g_2g_3=g_3g_2$ where $g_2$ is a 2-element  and $g_3$ is a 3-element of $G.$ It is easy to observe that $|\supp(g)|\leq |\supp(g_2)|+|\supp(g_3)|$. So either 
$|\supp(g_2)|>|\supp(D_2)|$ or $|\supp(g_3)|>|\supp(D_3)|$. In the former case
$g_2$ is not conjugate to an element  of $D_2$, so $\chi(g)=0$ for every \ir character $\chi$ in a 2-block with defect group $D_2,$ by Lemma \ref{jk1}. Similarly, consider the latter case. Here the result follows for $G=S_n$. If $G=A_n$ then the result follows from that for $S_n$ and the fact stated prior Lemma ref{ns7}. \dne

{\sc Comments:\, }  In view of Lemma \ref{23a}, it is desirable to determine the minimal defect group for every $n$ and $p=2$ or $3$.  If $p=2$ then the number $bp=2b$ must be of shape $k(k+1)/2$, so the minimal defect group of $S_n$  is a Sylow 2-subgroup of $S_{2b}$, where $k$ is the maximal integer such that $n-2b=k(k+1)/2$ for some $b>0$. If $p=3$ we only prove the existence of   a defect group $D$ with $|\supp(D)|\leq 2\sqrt{n}+4$ (Lemma \ref{3df}), which is sufficient for our purpose.

 \section{Minimal defect group of $S_n$ for $p=2$} 
 
The following lemma is obvious in view of the above comments. 

\bl{p2}  Let $p=2$ and $d\geq 0$ an integer. Then $d=|\supp(D)|$ for some defect group $D$ of a   $2$-block of $S_n$ \ii $d$ is even and $n-d=m(m+1)/2$ for some $m>0$.  In particular, 
 d is minimal \ii m is a maximal number such that $n-m(m+1)/2\geq 0$ is even. 
 \el

\bl{tr1}  Let $D$ be a minimal defect group of $S_n$ for $p=2$. If $n>13$ then  $|\supp(D)|< 3\sqrt{2n-20}$.
\el

\pf Let $T_m$ be the triangular diagram $[m,m-1\ld 1]$. Then $|T_{m+1}|-|T_m|=m+1$ as $T_{m+1}$   is obtained from $T_m$ by adding $m+1$ boxes. 

  Let $d=|\supp(D)|$ and let $m$ be the maximal integer such that $m(m+1)/2\leq n$. 
 Set  $a=n-m(m+1)/2$, so $a\leq m$. If $a$ is even then $d=a\leq m$. If $a=0$ then $d=0$, and the lemma is trivial.
 So we assume $a>0$. 
 

Suppose that  $a$ is odd. Then $T_m$ is not a 2-core of any diagram of  $S_n$.  So consider $T_{m-1}$.
Then $n-|T_{m-1}|=a+m$. If $m$ is odd then $a+m$ is even and then $d=a+m\leq 2m$.
Let $m$ even. Then $a\leq m-1$. Then $T_{m-1}$ is not a 2-core of any diagram of  $S_n$. Consider $T_{m-2}$.
Then $n-T_{m-2}=a+2m-1$ is even, so $d\leq 3m-2$.   

Therefore, $d\leq 3m-2$ in any case. 
We claim that $3m-2<3\sqrt{2n-20} $. Indeed, this is equivalent to $(3m-2)^2<18n-60$,
or $9 m^2-6m+64<18n$.  As $m(m+1)/2<n$, we have $9m^2+9m<18n$, it suffices to
show that $9 m^2-6m+64<9m^2+9m$, or $64<15m$. This is true if $m\geq 5$.
This holds if $n\geq 15$.  As $S_{15}$  has a block of defect 0, we are left with
$n=14$. If $n=14$  then $m=4$, $|T_m|=10$ and $d=4<3\sqrt{8}$, as claimed. \dne 

Recall that for $g\in S_n$ and a prime $p$ we denote by $g_p$ the element  such that $g=g_ph$, where $h\in \langle g \rangle$
and $|h|$ is not a multiple of $p$. One observes that $|\supp (g_p)|$
is the sum of cycle lengths divisible by $p$ in the cycle decomposition of $g$. 

\bl{tr2}  Let $G=S_n$ or $A_n$ where $n>13$. Let $R=\{g\in G: |\supp(g_2)|\geq  3\sqrt{2n-20}\}$.
Then there is a $2$-block B of G such that   $\chi(R)=0$ for every  \ir character $\chi$ of B. \el
 
\pf  Let $T$ be the triangle \di of maximal size $|T|\leq n$ such that $n-|T|$ is even. 

Suppose first that $|T|=n$.  Then, by Lemma \ref{d01}, $S_n$ has a character $\chi$ of 2-defect 0. Recall that $\chi$ is a unique character in a block it is contained in.  So for $G=S_n$  the statement  follows  from Lemma \ref{pr1}. If $G=A_n$ then $\chi|_{A_n}$ is the sum of two \ir characters of degree $\chi(1)/2$; moreover, each of them is of 2-defect 0 (indeed, an \ir character $\chi$ of $A_n$  is of 2-defect 0 \ii $\chi(1)$ is a multiple of $|S_n|_2$; as $|S_n|_2=2|A_n|_2$, it follows that $\chi(1)/2$  is a multiple of $|A_n|_2$). This implies the result for $A_n$.
 
 Let $|T|<n$ and let 
  $Y$ be any diagram of size $n$ containing $T$. Let $\chi$ be the \ir character of $S_n$ labeled by $Y$. Then $\chi$ belongs to a block $B$, say, whose defect group $D$ satisfies $|\supp(D)|=n-|T| $.
By Lemma \ref{tr1},  $|\supp(g_2)|\geq  3\sqrt{2n-20}>n-|T| $, so $g_2$ is not conjugate to an element of $D$.  By Lemma \ref{jk1}, $\chi(g)=0$ for every  \ir character $\chi$ of $ B$, whence  the result for $S_n$.
It is known that the defect groups of blocks of $A_n$ to which the \ir constituents of $\chi|_{A_n}$ belongs are $D\cap A_n$ 
and conjugate in $S_n$ (see comments prior Lemma  \ref{ns7}). So $g_2$ is not conjugate to an element of $D\cap A_n$, and the result follows as above for $S_n$. 

By Lemma \ref{ns7}, $Y$ can be chosen non-symmetric, so $\chi$ is \ir under restriction to $A_n$.   \dne

We say that the  element $g\in G\subseteq S_n$ has {\it maximal support} if $|\supp(g)|\geq |\supp(h)|$ whenever $h\in G$ and $|h|=|g|$.
One easily observes that if   $g$  is of maximal support in $G= A_n$ and $|g|$ is even  then $|\supp(g)|\geq n-3$.

\begin{corol}\label{cg6} Let $G=A_n$ or $S_n$ with $n>4$ and $n\neq 7,11,$ and let $M$ be the set of all $2$-elements of maximal support.   Then $\chi(M)=0$ for some \ir character $\chi$ of $G$. 
\end{corol}

\pf  Let $g\in G$ be a $2$-element  of maximal  support. Then   $|\supp(g)|\geq n-3$. If $n>13$ then $n-3\geq 3\sqrt{2n-20}$, and the result follows from Lemma \ref{tr2}. For $n\leq 13$ and  $n=5,6,8,10,12$ the result follows as   $A_n$ has a 2-block of defect 0. So we are left with $n=9,13$, which   can be  inspected by the character table of $G$  in~\cite{Atl}. 
Let $G=S_n$. Then we have to deal also with the cases with $n=5,8,12$. If $n=12$ the all elements of maximal support are in $A_{12}$. The cases with $n=5,8$ follows by inspection. 
  \dne

 Remark.  In $A_{7}$ the elements of order $2$ are not vanishing while all elements of order $4$ are vanishing. In $A_{11}$ the elements of order $2$ and maximal support are non-vanishing while all other 2-elements  are vanishing. In $A_{13}$ the elements 
 of maximal support and all 2-elements of order greater than 2 are vanishing, whereas all other 2-elements are non-vanishing
 (they form two conjugacy classes).

 For the use in Section  7 we compute the minimal numbers in the set $\{|\supp(D)|: D$ is a defect group of a 2-block of $S_n, 13<n<34\}$. Note that these are equal to $n-t$, where 
  $t$   is the maximal number of shape $m(m+1)/2$ such that $n-t$ is even.  
  
 Note that   $m(m+1)/2<34$ implies 
$m\leq 7$ so $t\in\{1,3,6,10,15,21,28\}$ for $n<34$. 

{\small
 \begin{center}
 TABLE 1
    \vspace{3pt}

   \noindent  \begin{tabular}{|c|c|c|c|c|c|c|c|c|c|c|c|c|c|c|c|c|c| }

        \hline
        $n$& $14$ & $16$ &$17$ &$18$  & $19$ & $20$ &$22$ &$23$& $24$ & $25$ &$26$ &$27$ &$29$&$30$ &$31$&$32$&$33$\\
        \hline
      $t$ &$10$ &  $10$ & $15$ &$10$&$15$& $10$&$10$&$21$&$10$&$21$&$10$&$21$  &$21$&$28$&$21$&$28$&$21$\\
\hline
        $n-t$ & $4$ &  $6$ &  $2$& $8$    & $4$ &  $10$ &  $12$& $2$ & $14$ &  $4$ &  $16$& $6$ & $8$ &  $2$  
& $10$& $4$& $12$
        \\
        \hline
  \end{tabular}
\end{center} }


  


\section{Defect groups for $p=3$ and vanishing elements}

We turn to the case $p=3$. 
 Consider the following diagrams, where $r,s\geq 0$ are integers and $r+s>0$: 

$C_{rs}=[r+2s,r+2s-2\ld r+2,r,r,r-1,r-1\ld 1,1]$ and  $C^-_{rs}=[r+1+2s,r+1+2(s-2)\ld r+2,r,r,r-1,r-1, \ld 1,1].$  
In particular,


$C_{0,s}=[2s,2s-2,2s-4\ld 4,2]$ ($s>0$) and $C^-_{0,s}=[2s+1,2s-1\ld 3,1]$  $(s\geq 0).$  

These diagrams appeared in \cite[p. 1159]{CCL}.  One can show that these exhaust all 3-cores in $S_n.$

Note that $|C_{rs}|=rs+r(r+1)+s(s+1)$ and $|C^-_{rs}|=(r+1)(s+1)+r(r+1)+s(s+1)$.
    In addition, $C_{rs},C_{sr}$ and $C^-_{rs},C^-_{sr}$ are transposes of each other. 

\medskip
Examples.  {\small $C_{0,4}=[8,6,4,2]$;  $~C_{0,3}^-= [7,5,3,1];$ $~C_{2,3}^-= [9,7,5,3,2,2,1,1]$;  $~C_{2,3}=[8,6,4,2,2,1,1]$}

\bigskip

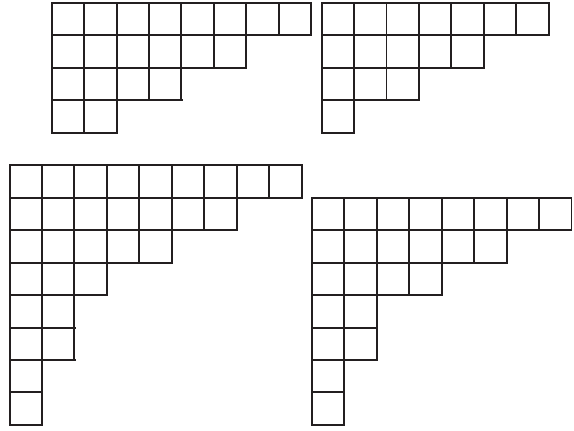
\centerline{  \yng(8,6,4,2) ~~~~~~~~~~~~~~~ 
\yng(7,5,3,1)} 

\bigskip
\centerline{
 \yng(9,7,5,3,2,2,1,1) ~~~~~~~~~~~~
 \yng(8,6,4,2,2,1,1)  }

 \bl{3df}   Let $G=S_n$, $n>13$. Then G has a $3$-block whose defect group support size is at most 
  $2\sqrt{n}+4$. 
 \el
 
\pf  We have to show that there is a 3-block with defect group $D$, say, such that $|\supp(D)|\leq 2\sqrt{n}+4$. We can assume that $S_n$ has no block of defect 0 (otherwise $D=1$ and the statement is trivial). 
 
Let $s$ be the maximal number such that $s^2\leq n$ (so $s\leq \sqrt{n}$). 
 Then  $n<(s+1)^2$. 
 Note that  $|C^-_{0,s-1}|=s^2$, $|C_{1,s-1}|=s^2+1$ and  $|C_{3,s-2}|=s^2+8$, and these three numbers have distinct residues modulo 3. 
 
 (i) Suppose that  $n> s^2+8$. As $n>13$, we have $s\geq 3$. Then we choose  the diagram  $C\in\{C^-_{0,s-2}$, $C_{1,s-1}$,   $C_{3,s-2}\}$ such that $3|(n-|C|)$. 
 Then there is a diagram $Y$ for $S_n$ such that 
 $|\tilde Y|=C$. 
 As $n<(s+1)^2$, we have  $n-|\tilde Y|\leq 2s\leq 2\sqrt{n-9}$.
 
 (ii)   Let $n\leq  s^2+8$, so $s>2$ (since $n>13$).  If the equality holds then $S_n$ has a block of defect 0. So $n\leq  s^2+7$.
 If $n=  s^2+7,s^2+4$  
 then we choose for $\tilde Y$ the \di $C_{1,s-1}$.  Then $n-|\tilde Y|
 \geq  s^2+7-(s^2+1)= 6$. 
 If $n=  s^2+6$ or $ s^2+3$ then  we choose for $\tilde Y$ the \di  $C_{1,s}$. Then $n-|\tilde Y|
 \geq s^2+6-s^2=6$. 
 
 Let $n=  s^2+5$ or $   s^2+ 2$.  Then  $n>13$ implies $s\geq 3$. we choose for $\tilde Y$  the 
  diagram  $C'\in\{C^-_{0,s-3}$, $C_{1,s-2}$,   $C_{3,s-3}\}$ such that 
   $3|(n-|C'|)$. Then $|\tilde Y|\geq (s-1)^2$, 
   and hence   $n-|\tilde Y|
 \geq s^2+5-(s-1)^2=2s+4\leq 2\sqrt{n}+4$. \dne
 
 
 \medskip
 {\sc Remark:\,} The bound  $2\sqrt{n}+4$ in Lemma \ref{3df} is not sharp.
 
 \medskip
 For  $g\in G= S_n$ denote by $g_3$ the element such that $g=g_3h=hg_3$, where $g_3$ is a $3$-element  and $|h|$ is coprime to $3$. 
 
  \begin{corol}\label{3dd}  Let $G= S_n$ or $A_n$, and 
  let $M=\{g\in G:|\supp(g_3)|> 2\sqrt{n}+4\}$.  
 Then there is an \ir character $\chi$ of G such that $\chi(M)=0$.

 \end{corol}

\pf By Lemma \ref{3df},  there is a 3-block of $S_n$ whose defect group support is at most $2\sqrt{n}+4$. So $g_3$ cannot lie in any defect group of this block. By Lemma \ref{jk1},  $\chi(g)=0$ for any \ir character in this block and any $g\in M$. (Recall that the defect groups of a block are conjugate, so their supports are of the same size.)  As the defect groups of 3-blocks of $A_n$ are the same as those of $S_n$, the result follows. \dne

 \medskip
 For further use  we record in Table 2
  the list of maximal 
  sizes of a core 
 for  $3$-blocks of $S_n$ for $n<51$ provided $S_n$ has no 3-block of defect 0. (In the latter case the size of some core equals $n$.) The numbers $n$ to be inspected are in the first row of the table, see Corollary \ref{i51} and the comment following it. 
In the  table $C$ is a core of maximal size $|C|\leq n$ with  $n\equiv |C|\pmod 3$, so $n-|C|=|\supp(D)|$,
where $D$ is a minimal defect group of $S_n$ for $p=3$.


{\small
 \begin{center}
 TABLE 2
    \vspace{3pt}

   \noindent  \begin{tabular}{|c|c|c|c|c|c|c|c|c|c|c|c|c|c|c|c|c|c|c| }

        \hline
        $n$& $7$ & $11$ &$13$ &$15$  & $18$ & $19$ &$23$ &$27$& $28$ & $29$ &$31$ &$35$ &$38$&$39$ &$43$&$45$&$47$&$48$\\
\hline
        $n-|C|$ & $3$ &  $3$ &  $3$& $3$    & $6$ &  $3$ &  $3$& $3$ & $3$ &  $3$ &  $6$& $3$ & $6$ &  $3$  
& $6$& $3$& $3$&$6$\\
\hline
      $|C|$ &$4$ &  $8$ & $10$ &$12$&$12$& $16$&$20$&$24$&$25$&$26$&$25$&$32$  &$32$&$36$&$37$&$42$&$44$&$42$
        \\
        \hline
  \end{tabular}
\end{center}
}
 
\medskip

 The table shows that for every $n<51$ there is an \ir character $\chi$ such that $\chi(g)=0$ whenever $g$ is a  3-element of order at least 9. Note that if $n=51$ then $|C|=42$ and    $n-|C|=9$. 

 \begin{lemma}\label{c51}  If $4<n<51$ then $S_n$ has a $3$-block of defect $0$ or $1$, unless $n=18,31,38,43,48$ where
there is a $3$-block of defect $2.$
\end{lemma}

\pf This immediately follows by inspection of Table 2. \dne

   \begin{propo}\label{ma3}  Let $ G=A_n$ with $n>4$ and  $n\neq 7,$  and let $e$ be the maximum order of a $3$-element of G.  
   Let $M $ be the set of all elements  of    order e. Then 
 $\chi(M)=0$ for some \ir character of $G$. 
 \end{propo}
 
\pf  Let $g\in M$. Observe that $|\supp(g)|> n/3$.  
 If $\frac{n}{3}\geq  2\sqrt{n}+2$ then the result follows from Lemma \ref{n44}. 
  Suppose that   $\frac{n}{3}<  2\sqrt{n}+2$. Then $n< 47$. 

If $27\leq n\leq 46$ then $e=27>16>2\sqrt{n}+2$. For $9\leq n\leq 26$ we have $e=9$ and   $|\supp(g)|> 9$.
Acording with Table 2, there is a 3-block of defect group support at most 6, so $g$ is not conjugate to any such  defect group.
Therefore, $\chi(g)=0$ for any \ir character in this block, whence the result. If $n=5,6,8$ then the groups $A_n,S_n$ have a 3-block
of defect 0, so we are left  with  $n=7$ as claimed.\dne

  \begin{propo}\label{3au}  Let $ G=S_n$ or $A_n$ with $n>4,$ $n\neq 7,$ and let $|G|=3^am$, where m is not a multiple of 3. Then there is a $3$-block of defect at most $ (a-1)/3$. 
 \end{propo}

\pf Suppose first that $\frac{n}{3}\geq |\supp(D)|,$ where $D$ is a   defect group of some 3-block $B$ of $G$.  Then
 $3|D|^3\leq |G|_3$ for $n\neq7$. Indeed, $G$ contains a 3-subgroup $X$, say, isomorphic to the wreath product of $D$ with $A_3$, and   $|X|=3|D|^3$. Let  $|D|=3^r$, so $r$ is the defect of the block   $B$. Then $|X|=3^{3r+1}$ and $3r+1\leq a$.
 This implies $r\leq (a-1)/3$.

Suppose that 
$n\geq 58$. Then 
$\frac{n}{3}>  2\sqrt{n}+4$. By Lemma~\ref{3df} \,there is a 3-block of $G$ whose defect group satisfies 
 $|\supp(D)|\leq2\sqrt{n}+4$. So the result follows.  Let $n<  58$. By Table 2, $\frac{n}{3}\geq   |\supp(D)|$ unless $n=7$.\dne

  



 \subsection{Elements of maximal support whose order is divisible by 3}


  Let $G= S_n$ or $A_n$ and $g\in G$. 
 Recall that $g$ is of {\it maximal support} in $G$ if ${|\rm supp}(g)|\geq |{\rm supp}(h)|$ for any element $h$ of the same order as $g.$ In other words, the  number of fixed points of such an $h$  does not exceed the number of fixed points of $g$.  
 
 Note that for $g$ in $G$ to be of maximal  support depends on whether $G=S_n$ or $A_n$. Say, if $g\in A_6$ is an involution (a double transposition) then $g$ is of maximal  support in $A_6$ but not in $S_6$. The following is easily verified: 
 
\bl{ms1}  Let $G=  S_n$ or $A_n$ and $g\in G$.  Suppose that $g$ is of maximal  support. \\[7pt]
{\rm (1):}  Let $G=  S_n$ and  let $|g|$ be even. Then $|{\rm supp}(g)|\geq n-1$. \\[7pt]
{\rm (2):} Let   $g\in A_n$ and let $|g|$ be even. Then $|{\rm supp}(g)|\geq n- 3$. \\[7pt]
{\rm (3):} Suppose that  $|g| $ is a multiple of $3$. Then $|{\rm supp}(g)|\geq n- 2$.  \\[7pt]
{\rm (4):} Suppose that  $|g| $ is a multiple of $6$. If $G=A_n$ then $|{\rm supp}(g)|\geq n- 2$, and if $G=S_n$ then 
$|{\rm supp}(g)|\geq n- 1$. \el

  \bl{23b}  Let $G=S_n$ or  $ A_n$ with $n>13$ 
  and let $g=g_2g_3=g_3g_2\in G$, where $|g_2|=2^\al$  and $|g_3|=3^\beta$. Suppose that $|\supp(g_2)|>  3\sqrt{2n-20}$ or $|\supp(g_3)|>  2\sqrt{n}+4$. Then there is an \ir character $\chi$ of $G$ such that $\chi(g)=0$. 
\el

\pf If  $|\supp(g_2)|>  3\sqrt{2n-20}$ then $\chi(g)=0$ for some character $\chi$ of $S_n$ by Lemma \ref{tr2}.  
If  $|\supp(g_3)|>  2\sqrt{n}+4$ then $\chi(g)=0$ 
for some \ir character of   $G$ by Corollary \ref{3dd}. 
So the lemma follows for $G=S_n$. 

Let $G=A_n$. If $\chi|_{A_n}$ is \ir then we are done.  Otherwise, $\chi|_{A_n}$ is the sum of two \ir characters of equal degree.  If $D=1$ then  $\chi$ is of defect 0, which is equivalent to saying that $\chi(1)$ is a multiple of $|S_n|_3=|A_n|_3$ (see \cite[Theorem 86.3]{CR}). 
Suppose that $D\neq 1$, that is, the block $\chi$ belongs to is not of defect 0.
Then, by Lemma \ref{ns7}, this block has an \ir character $\chi'$ labeled by a non-symmetric Young diagram which is  therefore \ir under restriction to $A_n$. \dne

 Obviously, Theorem \ref{main2} follows from the following result:

 \begin{theo}\label{mt2}  Let $G=S_n$ or $A_n$ with $n>4$ and let $g\in G.$ 
 Suppose that g is of maximal  support in G. Then $\chi(g)=0$ for some \ir character of G, unless 
$n=7$ or $11$.   
\end{theo}

{\it Proof.} In view of Theorem \ref{tc1}, we can assume that $|g|=2^\al 3^\beta,$ and by Corollary \ref{cg6}, that $\beta\neq 0$.
(Here $\al,\beta\geq 0$ are integers.)
By Lemma \ref{sd9}, we can assume that $g\in A_n$. Let $g=g_2g_3=g_3g_2$, where $|g_2|=2^\al$, $|g_3|=3^\beta$.
By Lemma \ref{jk1}, we are done if  $|\supp(g_2)|>|\supp(D_2)|$ or $|\supp(g_3)|>|\supp(D_3)|$ for some defect groups
$D_2,D_3$  of some 2-block and 3-block of $G$, respectively.  So we choose $D_2,D_3$ to be minimal.

Suppose first that $n>13$. By Lemma \ref{23b}, we can assume that $|\supp(g_2)|\leq  3\sqrt{2n-20}$ and $|\supp(g_3)|\leq  2\sqrt{n}+4$. As $g$ is  of maximal  support, we have $n-2\leq |\supp(g)|\leq |\supp(g_2)|+|\supp(g_3)|$ (Lemma \ref{k9}), whence
$n-2\leq 3\sqrt{2n-20}+2\sqrt{n}+4  $.  This implies   $n<45$.  


Let $n<45$.
One observes in Table 2   that $|\supp(D_3)|\leq 6$ for $n$ with $4<n<45$. 
 (If $D_3=1$ then $|\supp(D_3)|=0$.) By Lemma \ref{jk1}, we are done unless $|\supp (g_3)|\leq 6$. If so,
 we have $n-2\leq |\supp(g)|\leq |\supp(g_2)|+6$, whence
 $|\supp(g_2)|\geq n-8$. By Lemma \ref{tr2}, we are done if $|\supp(g_2)|\geq   3\sqrt{2n-20}$. The latter  implies 
  $ 3\sqrt{2n-20}> n-8$, whence  
 $n\leq 23$.  
 The cases $13<n\leq 23$ are verified by inspection, and we ignore the values of $n$ such that   $A_n$ and $S_n$ has a 3-block of defect 0. For the remaining cases $n=23,19,18,15$ we have 
 $|\supp(D_2)| =2, 4,8,0$, respectively, for $S_n$. As the defect groups of 2-blocks for $A_n$ are $D_2\cap A_n$ when $D_2$ runs over the defect groups of 2-blocks of $S_n$, we have $|\supp(D_2)|\leq 8$, so the result follows.

The cases with $n< 13$ can be examined from their character tables \cite{Atl}.   If $14>n>4$ and $n\neq 7,11$ then $A_n$ has a block of 3-defect 0. The non-vanishing elements     $g\in G=A_{11}$  are in class $2B$ and $3A$ in notation of \cite{Atl}, and those in $3A$  are not of maximal support.  Let $G=A_7$. Then an element of  maximal support is in class $3B$ or $ 6A$, the latter  is  
non-vanishing, see ~\cite{Atl}. If $n=13$ then the character table is available in the GAP library. One observes that the only non-vanishing elements in $A_{13}$
are in classes $2A,2B,3A$, and hence not of maximal support. This completes the proof. 
\dne



\medskip
Denote by $\omega(G)$ the set of element orders in the group $G.$

\bl{e6}  Let $g\in G=A_n$.  Suppose that $2|g|,3|g|\notin \omega(G)$. \\[5pt]
$(1)$ We have $|\supp(g)|\geq n/2;$\\[5pt]
$(2)$ Let $g=g_2g_3=g_3g_2$, where $|g_2|=2^\al$, $|g_3|=3^\beta$, $\al,\beta>0$.
Then $n< \min\{2^{\al+1}+|\supp (g_3)|+2, |\supp(g_2)|+3^{\beta+1}\}$. 
\el

\pf  (1) Suppose the contrary. Then $|\supp(g)|\leq \frac{n}{2}-1$. Let $c_1\ld c_k,1\ld 1 $ be the cycle lengths of $g$, where $c_1\ld c_k>1$. So $|\supp(g)|=c_1+\cdots +c_k$. Let $h$ be the permutation 
with cycle lengths $2c_1\ld 2c_k,1\ld 1 $. Then $|h|=2|g|$
and $|\supp(h)|=2|\supp(g)|\leq n-2$.
If $h\notin A_n$ then  replace $h$ by the permutation $h'$   obtained from $h$ by adding  a 2-cycle. 
Then $h'\in A_n$ and $|h'|=2|g|$,
which is a contradiction. 

(2) Suppose the contrary. If  $n\geq 2^{\al+1}+|\supp (g_3)|+2$ then consider an element $g_2'g_3$, where $|g_2'|=2^{\al +1}$ and  $| \supp(g_2')\cap \supp(g_3)|=0$. Such element $g_2'\in A_n$ exists as $n- |\supp (g_3)|\geq  2^{\al+1}+2$.  Then $|g_2'g_3|=2|g|$, a contradiction. Similarly, if $n\geq 3^{\beta+1}+|\supp (g_2)|$ then there is an element $g_3'$ such that $ |g_3'|=3^{\beta+1}$ and $| \supp(g_2)\cap \supp(g_3')|=0$. Then $A_n$ contains an  element $g_2g_3'$ of order $3|g|$.  
 \dne

{\it Proof of Theorem} \ref{pp7}.
  Let $|g|=2^\al 3^\beta $. Suppose that $n>13$. By Lemma \ref{23b}, we can assume that 
$|\supp(g_2)|< 3\sqrt{2n-20}$
and  $|\supp(g_3)|< 2\sqrt{n}+4$. By Lemma \ref{e6}, $|\supp(g)|\geq n/2$ and hence
either $|\supp(g_2)|\geq n/4$ or $|\supp(g_3)|\geq n/4$. So 
either $3\sqrt{2n-20}> n/4$ or $2\sqrt{n}+4> n/4$. (Note that $3\sqrt{2n-20}>  2\sqrt{n}+4$ for $n>15 $.)

Suppose  that $n>15$ and $3\sqrt{2n-20}> n/4$. Then 
$n\leq 297$.  Then  $ |\supp(g_2)|< 3\sqrt{2n-20}<72$ for $n\leq 297$ and $|\supp(g_3)|< 2\sqrt{n}+4<39$. Therefore, $2^\al\leq 64,3^\beta\leq 27$.   Note that              $ |\supp(g)|= |\supp(g_2)|+ |\supp(g_3)|-| \supp(g_2)\cap \supp(g_3)|$. Therefore,  $ |\supp(g)|\leq |\supp(g_2)|+ |\supp(g_3)|$. 
By Lemma \ref{e6}(2), we have $n\leq  3^4+72=153$. 

This argument can be repeated: $ |\supp(g_2)|< 3\sqrt{2n-20}<51$ for $n\leq 153$ and $|\supp(g_3)|< 2\sqrt{n}+4<29$. Again by  Lemma \ref{e6}(2), we have $n\leq 3^4+50=132$. 
 
  Again $ |\supp(g_2)|< 3\sqrt{2n-20}<47$ for $n\leq 132$ and $|\supp(g_3)|< 2\sqrt{n}+4<27$, whence
 $3^\beta\leq 9$.  Then $n\leq 3^3+46=73$. 

 Once more  $ |\supp(g_2)|< 3\sqrt{2n-20}<34$ for $n\leq 73$ and $|\supp(g_3)|< 2\sqrt{n}+4<22$, whence  $n\leq
 3^3+32=59$. 
 
  Finally,    $ |\supp(g_2)|< 3\sqrt{2n-20}<30$ for $n\leq 59$ and $|\supp(g_3)|< 2\sqrt{n}+4<20$. Therefore, $2^\al\leq 16$ and  $|\supp(g_3)| \leq 18$. By  Lemma \ref{e6}(2), we have $n\leq 2^5+18+2=52$.
  
  If $n=52$ then $G$ has a 3-block of defect 0, if $n=51$ then there is a 3-block of $G$ with defect group support 9 (see comment after Table 2). In view of Lemma \ref{jk1}, we can assume that $|\supp(g_3)|\leq 9$, and then $n\leq 2^5+2+9=43$ by Lemma \ref{e6}(2). For $n\leq 43$ the group $G=A_n$ has a 3-block $B$ with defect group support at most 6
  (Table 2). Then, as above,  $|\supp(g_3)|\leq 6$, so  $|g_3|=3$, $\beta=1$.
    Then Lemma \ref{e6}(2) yields $n\leq 2^5+2+3=37$.   
    
    For $n\leq 37$ one observes from Table 2 that $G$ has a 2-block with  defect group support at most 16. Then
    $|\supp(g_2)|\leq 16$, and 
      $n\leq 3^2+16\leq 25$ by  Lemma \ref{e6}(2). 
     For $n=25,24$ the group $G$ has a 3-block of defect 0, so these values are ruled out too. If $n=23$ then 
    $G$ has a 2-block with  defect group support 2, whence $|\supp(g_2)|\leq 2$,  which contradicts the inequality $n\leq 9+4$ in Lemma   \ref{e6}(2).
     For $n=22,21,20$ the group $G$ has a 3-block of defect 0. Let $n=19$. Then 
    $G$ has a 2-block with  defect group support 4, whence $|\supp(g_2)|\leq 4$,  which contradicts the inequality $n\leq 9+6$ in Lemma   \ref{e6}(2).   So $n\leq 18$.
  
  Let $n=18$. Suppose that $|g_2|=8$. Consider the diagram $Y=[12,2,2,1,1]$:
  
  \medskip
  \centerline{  \yng(12,2,2,1,1)} 
   Observe that there is a unique way to delete an $8$-rim from $Y$. Let $\chi_Y$ be the \ir character of $S_{18}$. Applying the Murnahgan-Nakayama formula (\ref{eq MN}) with $c=8$ we   get
  $\chi_Y(g)=\chi_{Y_1}(h)$, where $Y_1=Y\setminus \nu=[4,2,2,1,1]$, $h\in S_{10}$  is obtained from $g$ by removing the $8$-cycle   and $\chi_{Y_1}$ is an \ir character of $S_{10}$. 
  There is no way to delete a $3$-rim from $Y_1$; so $\chi_{Y_1}$ is of 3-defect 0.  As $|h|$ is a multiple of 3, we have $\chi_{Y_1}(h)=0$, and hence $\chi_Y(g)=0$. 
  
 So the result follows if  $|g_2|=8$. Let $|g_2|\leq 4$. As $|\supp(g_3)|\leq 6$, by Lemma \ref{e6}(2) we have $n\leq 6+2^3+2=16$, a contradiction.  
 
 If  $n=17,16,14$ then $G$ has a 3-block of defect 0, and if $n=15,12$ then $G$ has a 2-block of defect 0.
 If   $n\leq 13$ then the result follows by inspection of the character table of $G$.  \dne


\begin{thebibliography}{hhhhhh}
 
\bibitem{CCL} X. Chen, J. Cossey, M.L. Lewis and H.P. Tong-Viet, Blocks of small defect in alternating groups and squares of Brauer character degrees, J. Group Theory 20   (2017), 1155 - 1173. 
 
 \bibitem{Atl} 
  J. Conway, R. Curtis, S. Norton, R. Parker and R. Wilson,
  \emph{Atlas of Finite Groups}. Clarendon Press, Oxford, 1985.

 
 \bibitem{CR} Ch. Curtis and I. Reiner, {\it Representation theory of finite groups and associative algebras},  Wiley\&Son, New York, 1962.
 
\bibitem{DN} S.  Dolfi, G.  Navarro, E. Pacifici, L. Sanus and  Pham Huu Tiep,
Non-vanishing elements of finite groups. 
J. Algebra  323  (2010),  no. 2, 540 -– 545. 

 
 
 
 \bibitem{GO} A. Granville and K. Ono, Defect zero $p$-blocks for finite simple groups, 
 Trans. Amer. Math. Soc. 348(1996),  331 –- 347.
 
  \bibitem{GutBo} I. Gutman and B. Borovicanin, Nullity of Graphs: An updated Survey, \\ http://www.mi.sanu.ac.rs/projects/ZbR14-22.pdf
 

\bibitem{Is1} M. Isaacs,  {\it Character theory of finite groups}, Academic Press, N.Y. 1978.
 
\bibitem{INW} M. Isaacs, G. Navarro and  T. Wolf, Finite group elements where no \ir character vanishes, J. Algebra
222(1999), 413 -- 423.


\bibitem{JK} G. James and Kerber, {\it  The representation  theory of the symmetric group}, Addison-Wesley, London, 1981.


\bibitem{K} A. Klyachko, Modular forms and representations of symmetric groups, J. Soviet Math. 26(1984), 1879 - 1887.


\bibitem{Lo} L. Lov\'asz,  {\it Spectra of graphs with transitive  groups}, Period. Mat. Hungar. 6 (1975) 191-196. 



\bibitem{MNO} G. Malle,  G. Navarro and J. Olsson, Zeros of characters of finite groups,  J. Group Theory, 4 (2000), 353-368.

\bibitem{na1} G. Navarro, {\it  Characters and blocks of finite groups}, London Mathematical Society Lecture Notes, Cambridge Univ. Press, Cambridge, 1998

\bibitem{AStreit}
 Andrew Streitwieser,  {\it Molecular Orbital Theory for Organic Chemists}, Wiley, New York, 1961.
 


\end{thebibliography}
\end{document}